\documentclass[11pt]{article}
\usepackage{amsmath, latexsym, amssymb,pstricks}
\usepackage[dvipdfm]{graphicx}

\makeatletter
\newfont{\footsc}{cmcsc10 at 8truept}
\newfont{\footbf}{cmbx10 at 8truept}
\newfont{\footrm}{cmr10 at 10truept}
\renewcommand{\ps@plain}
\makeatother 

\numberwithin{equation}{section}
\newtheorem{theorem}{Theorem}[section]
\newtheorem{thm}[theorem]{Theorem}

\newtheorem{rem}[theorem]{Remark}

\newenvironment{proof}{\medskip\noindent{\it Proof.\ }}{\hfill \mbox{$\Box$}\medskip}

\begin{document}


\title{Maximizing the number of maximal independent \\ sets  of a fixed size}

\author{Chunwei Song\footnote{School of Mathematical Sciences \& LMAM, Peking University, Beijing 100871, P.R. China {\tt csong@math.pku.edu.cn},
{\tt byao@pku.edu.cn}}, Bowen Yao\footnotemark[2] \thanks{The authors were partially supported by NSF of China grant \#11771246.}}

\date{}

\maketitle
\medskip

\small \emph{Mathematics Subject Classification 2020:} 05C35, 05C69, 05D99, 05C31

\emph{Keywords: } maximal independent set,
extremal graph, disjoint union,
Tur\'an graph, maximal cliques, maximal independence polynomial

\section*{Abstract}
For a fixed graph $G$, a maximal independent set is an independent set
that is not a proper subset of any other independent set.
P. Erd{\H o}s, and independently, J. W. Moon and L. Moser (\emph{Israel J. Math.}, 3 (1965): 23–28),
and R. E. Miller and D. E. Muller (\emph{IBM Res. Rep.}, (1960): RC-240),
determined the maximum number of maximal independent sets in a graph on $n$ vertices,
as well as the extremal graphs.
In this paper we maximize the number of maximal independent sets
of a fixed size for all graphs of order $n$ and determine the extremal graphs.
Our result generalizes the classical result.


\section{Introduction}
Throughout this paper, we consider finite simple connected graphs.
Let $G= (V(G), E(G))=(V, E)$ be such a graph with vertex set $V$ and edge set $E$.
Below are some graph theory concepts and notation needed in this paper.
Readers are suggested to refer to \cite{Bol98}, \cite{Die00} or \cite{Wes96} for  terminologies
not specified here.

An independent set (or stable set) of $G$ is a set of pairwise nonadjacent vertices.
In order for a set of vertices $U \subseteq V$ to be a \emph{maximal independent set} (abbr. \emph{MIS}) of $G$,
we require the set $U$ to (1) be independent, and (2) have no \emph{strictly super} independent set $W$ such that
$U \subseteq W \subseteq V$. The set of neighbors of $v \in V$ is denoted by $N(v)$, or if necessary by $N_G(v)$.
$G'=(V',E')$ is an \emph{induced subgraph} of $G$ if $V' \subseteq V$ and $E'$ consists of all edges
$uv \in E$ with $u, v \in V'$. We usually write $G':=G[V']$.

A $k$-partite graph is a graph whose graph vertices can be partitioned into $k$
disjoint sets so that no two vertices within the same set are adjacent.
A $k$-partite graph is said to be complete if every pair of graph vertices in
the $k$ sets are adjacent.
The \emph{Tur\'an graph} $T_{n,k}$ is the complete
$k$-partite graph with $n$ vertices whose partite sets differ in size by at most one \cite{Bol98}.
We also define $G_1+G_2$ as the graph consisting of the disjoint union of two graphs
$G_1$ and $G_2$.  The Tur\'an graph $T_{n,k}$ and the operation of disjoint union will
be useful in the extremal graph structure described in our main theorem (Theorem \ref{thm:general}).

Given $G$, let $i_t(G)$ be the number of independent sets of size $t$
in $G$ and let $i(G)=\sum_{t \geq 0} i_t(G)$ be the total number of independent sets.
While there have been many extremal results on $i(G)$ and $i_t(G)$ over various families of graphs
(see e.g. \cite{Gal11,GLS,GZ11,SSZ19,Zha10}),
it makes sense to investigate parallel theories on the MIS's,
the independent sets that are not covered by bigger ones.
Let $i^{max}_t(G)$ be the number of maximal independent sets of size $t$
in $G$ and let $i^{max}(G)=\sum_{t \geq 0} i^{max}_t(G)$ be the total number of maximal independent sets.
For arbitrary graphs $G$ on $n$ vertices,
P. Erd{\H o}s (see \cite{Fur87}), and independently, Moon and  Moser  \cite{MM65},
and  Miller and  Muller \cite{MM60}
determined  $i^{max}(G)$ as well as the extremal graphs.


Nonetheless, studies on the number of maximal independent sets seem to be less adequate
(see for instance \cite{Fur87, Lin18, LJ19, Wil86, Woo11}).

The \emph{maximal independence polynomial} is defined by \cite{HMS17}
$$I_\mathrm{max} (G; x) := \sum_{U:  \ U \ \text{is an MIS of} \ G} x^{|U|}. $$
By definition, $i^{max}_t(G)=[x^t] I_\mathrm{max} (G; x)$, where
by usual convention $[x^k] f(x)$ represents the coefficient of
$x^k$ in the polynomial or series $f(x)$.

In this note we maximize $i^{max}_t(G)$ graphs $G$ on $n$ vertices.

\begin{thm}\label{thm:general}
  Assume $n=qt+r$, where $0 \leq r <t$.
  For all graphs $G$ on $n$ vertices,
  we have
  \begin{eqnarray}
  & i^{max}_t(G) \leq q^{t-r} (q+1)^{r}. \label{formula:general}
   \end{eqnarray}
  Furthermore, let $H=(t-r) K_q + r K_{q+1}$, i.e. disjoint union of $t$
  cliques of the specified orders, then $H$ is the unique extremal graph.
\end{thm}

\begin{rem}
  If $n$ is a multiple of 3, \eqref{formula:general} shows that
  $i^{max}_{\frac{n}{3}}(G) \leq 3^{\frac{n}{3}}$
  and it implies the main theorem in \cite{MM65} which says that
  $i^{max}(G) \leq 3^{\frac{n}{3}}$ when 3 divides $n$.
  If $n=3k-1$,  \eqref{formula:general} gives that
  $i^{max}_{k}(G) \leq 2 \cdot 3^{k-1}$.
  If $n=3k+1$,  \eqref{formula:general} gives that
  $i^{max}_{k}(G) \leq 2^2 \cdot 3^{k-1}$
  and  $i^{max}_{k+1}(G) \leq 2^2 \cdot 3^{k-1}$.
  Each case above strengthens a respective case of \cite[Theorem 1]{MM65}.
  (Note that the celebrated result of \cite[Theorem 1]{MM65} says that
  each extremal graph actually has maximal independent sets of only one certain size.)
\end{rem}


\section{Proof}

For convenience we work with the complementary graph, and count cliques instead of independent sets.
That is, we show that for all graphs $G$ on $n=qt+r$ vertices, where $0 \leq r <t$,
the number of $t$-maximal cliques in $G$ is no more than $f(n,t):=q^{t-r} (q+1)^{r}$.
Furthermore, the Tur\'{a}n graph $T_{n,t}=K_{q,\cdots,q,q+1,\cdots, q+1}$
is the unique extremal graph. We achieve this by induction on $n+t$.
Keep in mind that the proposed extremal value $f(n,t)=q^{t-r} (q+1)^{r}$
strictly increases with $n+t$.

\begin{proof}
Since the cases that $n<t$ or $t=1$ are trivial,
without loss of generality, we assume that $n \geq t \geq 2$.

{\bf{Case 1.}} $r>0$.
We start with the case that is more convenient to phrase and {Case 2} will be easier to understand.

{\bf{Subcase 1a}}. $r>0$ and $\delta(G)\geq n-q$. Let $v_1, v_2, \cdots, v_t$ be arbitrarily selected.
Note that
$$|V - \bigcap_{i=1}^t N(v_i)|=|\bigcup_{i=1}^t (V - N(v_i)) | \leq \sum_{i=1}^t | V - N(v_i)| \leq qt < qt+r = n = |V|.$$
Hence $\bigcap_{i=1}^tN(v_i)\neq\varnothing$, i.e., every $t$ vertices in $G$ has a common neighbor.
 Thus any maximal clique is larger than $K_t$, so that  $G$ has no $t$-maximal cliques,
implying that this subcase needs not be considered in order to maximize $i^{max}_t(G)$.

{\bf{Subcase 1b}}.  $r>0$ and $\delta(G)\leq n-q-1$.
Choose a vertex $v$ such that $d(v)=\delta(G) \leq n-q-1$.

Let $\mathcal{A}$ be the set of $t$-maximal cliques in $G$ which contains $v$, and $\mathcal{B}$ be the set of $t$-maximal cliques in $G$ which does not contain $v$.

Every $t$-maximal clique in $\mathcal{B}$ is a $t$-maximal clique of $G-\{v\}$.
By induction hypothesis, as $n-1=qt+r-1$, $|\mathcal{B}|\leq q^{t-r+1}(q+1)^{r-1}$.

To calculate $|\mathcal{A}|$, consider the subgraph of $G$ induced by the neighbors of $v$.
Every $t$-maximal clique in $\mathcal{A}$ corresponds to a $(t-1)$-maximal clique of $G[N(v)]$.
As $|V(G[N(v)])| \leq n-q-1=qt+r-q-1=q(t-1)+r-1$, inductively,
$|\mathcal{A}|\leq q^{(t-1)-(r-1)}(q+1)^{r-1}=q^{t-r}(q+1)^{r-1}$.

Therefore, the total number of $t$-maximal cliques in $G$ is bounded by
$$|\mathcal{B}|+|\mathcal{A}| \leq q^{t-r+1}(q+1)^{r-1} + q^{t-r}(q+1)^{r-1} = q^{t-r}(q+1)^r=f(n,t).$$

The above equality holds if and only if the following conditions are simultaneously met.

i). $G-v$ is a $t$-cliques extremal graph of order $n-1$. Inductively, this requires $G-v=T_{n-1,t}=K_{q,\cdots,q,q+1,\cdots, q+1}$
with $t-r+1$ partite sets of size $q$ and $r-1$ partite sets of size $q+1$.

ii). $d(v)=n-q-1=qt+r-q-1=q(t-r)+(q+1)(r-1)$ and $G[N(v)]$ is a $(t-1)$-cliques extremal graph of order $n-q-1$. This requires that
$G[N(v)]=T_{n-q-1,t-1}=K_{q,\cdots,q,q+1,\cdots, q+1}$ with $t-r$ partite sets of size $q$ and $r-1$ partite sets of size $q+1$.

Putting i) and ii) together, clearly, the neighbors of $v$ are precisely $t-r$ of the total $t-r+1$
partite sets of size $q$ and $r-1$ partite sets of size $q+1$ in $G-v$. Thus $G=T_{n,t}$ is the unique extremal
graph in this case.

{\bf{Case 2.}} $r=0$, so that $n=qt$.
The case $r=0$ is similar to Case 1, with only slight differences in the calculation.

{\bf{Subcase 2a}}. $r=0$ and $\delta(G)\geq n-q+1$. For any $t$ vertices $v_1,v_2,\cdots,v_t$, as
$$|V - \bigcap_{i=1}^t N(v_i)| =|\bigcup_{i=1}^t (V - N(v_i)) \leq (q-1)t < |V|,$$
they must have a common neighbor. Thus $G$ has no $t$-maximal cliques.

{\bf{Subcase 2b}}.  $r=0$ and $\delta(G)\leq n-q$. Choose a vertex $v$ such that $d(v)=\delta(G)\leq n-q$.

Define $\mathcal{A}$ and $\mathcal{B}$ as in Case 1.
As $n-1=qt-1=(q-1)t+t-1$, by similar arguments, inductively,
$|\mathcal{B}| \leq (q-1) q^{t-1}$.

On the other hand, as $|V(G[N(v)])| \leq n-q=q(t-1)$, by induction,
we have $|\mathcal{A}|\leq q^{t-1}$.

Thus the total number of $t$-maximal cliques in $G$ is limited by
$$|\mathcal{B}|+|\mathcal{A}|\leq (q-1) q^{t-1} + q^{t-1}=q^t=f(n,t).$$

The above extremal value is achieved if and only if the following conditions are simultaneously met.

i). $G-v$ is a $t$-cliques extremal graph of order $n-1$. This means $G-v=T_{n-1,t}=K_{q,\cdots,q,q-1}$
with $t-1$ partite sets of size $q$ and $1$ partite set of size $q-1$.

ii). $d(v)=n-q=q(t-1)$ and $G[N(v)]$ is a $(t-1)$-cliques extremal graph of order $n-q$. This requires that
$G[N(v)]=T_{n-q,t-1}=K_{q,\cdots,q}$ with $t-1$ partite sets of size $q$.

Altogether, it is implied that $G=T_{n,t}$ is the unique extremal
graph in Case 2 as well.

\end{proof}

\begin{rem}
Theorem \ref{thm:general} says that the Tur\'{a}n graph $T_{n,t}=K_{q,\cdots,q,q+1,\cdots, q+1}$
is the unique extremal graph of the complementary scenario.  That is, $T_{n,t}$
has the maximum number of maximal cliques of size $t$.
Equivalently,
$\overline{T_{n,t}} = H=(t-r) K_q + r K_{q+1}$, disjoint union of $t$
cliques of most possibly balanced sizes,  is the unique extremal graph
that has the maximum number of maximal independent sets of size $t$.
\end{rem}



\end{document}